%
%

\input amstex
\input epsf
\documentstyle{amsppt}
\topmatter 
\NoBlackBoxes
\define\mathcal{\Cal}
\define\disj{\coprod}

\define\ind{\operatorname{ind}}
\define\zq{\Bbb{Z}Q^{op}}

\define\Ext{\operatorname{Ext}}
\define\Hom{\operatorname{Hom}}
\define\Dee{\Cal{D}^b}
\define\Cee{\Cal{C}}
\define\el{\Cal{L}}

\define\mathbb{\Bbb}
\define\aphi{\Phi_{\geq -1}}
\define\mphi{\Phi_{\geq -1}^m}

\define\st{\Cal{T}}

\define \qo{Q^{op}}
\define \qb{Q_{bip}}

\define\zqb{\mathbb{Z}\qb^{op}}
\define\cm{{\Cee_m(\Phi)}}

\define\mpsi{\Psi_{\geq -1}^m}
\define\dd{{\Cal{D}^b(\Phi)}}
\define\dpsi{{\Cal{D}^b(\Psi)}}
\define\ceem{\cm}

\redefine\ss{\Cal{S}}

\title Defining an $m$-cluster category 
\endtitle

\rightheadtext{Defining an $m$-cluster category}
\author Hugh Thomas\endauthor

\date September 2007 \enddate


\address Department of Mathematics and Statistics, University of 
New Brunswick, Fredericton NB, E3B 5A3, Canada 
\endaddress
\email hugh\@erdos.math.unb.ca \endemail

\abstract
We show that a certain orbit category considered by Keller 
encodes the combinatorics of
the $m$-clusters of Fomin and Reading in a fashion similar to the
way the cluster category of Buan, Marsh, Reineke, Reiten, and Todorov 
encodes the combinatorics of the clusters of Fomin and Zelevinsky.  
This allows us to give type-uniform proofs of certain results
of Fomin and Reading in the simply laced cases.  
\endabstract
\endtopmatter

\document

For $\Phi$ any root system, Fomin and Zelevinsky [FZ] 
define a {\it cluster complex}
$\Delta(\Phi)$,
a simplicial complex on $\aphi$, the {\it almost positive roots} of $\Phi$.  
Its facets (maximal faces) are called {\it clusters}.
In [BM+], starting in the more general context of 
a finite dimensional hereditary algebra $H$ over
a field $K$, 
Buan et al.~define a {\it cluster category} $\Cee(H)=
\Dee(H)/\tau^{-1}[1]$.  
($\Dee(H)$ is the bounded derived category of representations of $H$; 
more will be said below about it, its shift functor [1], and
its Auslander-Reiten translate $\tau$.)
The cluster category $\Cee(H)$ is a triangulated Krull-Schmidt category.
We will be mainly interested in the case where $H$
is a path algebra associated to the simply laced root system $\Phi$,
in which case we write $\Cee(\Phi)$ for $\Cee(H)$.  
There is a bijection $V$ taking $\aphi$ to the indecomposables of 
$\Cee(\Phi)$.  
A {\it (cluster)-tilting set} in $\Cee(\Phi)$ is a maximal set $\ss$ of indecomposables
such that $\Ext^1_{\Cee(\Phi)}(X,Y)=0$ for all $X,Y\in \ss$.  
$\Cee(\Phi)$ encodes the combinatorics of $\Delta(\Phi)$ in the sense that 
the clusters of $\Phi$ correspond bijectively to the tilting sets of 
$\Cee(\Phi)$ under the map $V$.  

Tilting sets in $\Cee(\Phi)$ always have cardinality $n$, the
rank of $\Phi$.  An {\it almost complete tilting set} is a set $\st$ 
of $n-1$ indecomposables such that $\Ext^1_{\Cee(\Phi)}(X,Y)=0$ 
for $X,Y\in \st$.  
A {\it complement} for $\st$ is an indecomposable $M$ such that $\st\cup\{M\}$
is a tilting set.  A tilting set always has exactly two complements.
(This was shown from the cluster perspective in [FZ] and from the
representation theoretic perspective in [BM+].) 

In [FR], Fomin and
Reading introduced a generalization of clusters known as
$m$-clusters, for $m\in \mathbb{N}$.
When $m=1$, the classical clusters
are recovered.  The $m$-cluster complex $\Delta_m(\Phi)$ is a simplicial
complex on a set of {\it coloured roots} $\mphi$.  It has been studied
further in [AT1, T, AT2].
The facets of $\Delta_m(\Phi)$
 are known as $m$-clusters.  
The goal of this paper is to show that the category 
$\Cee_m(\Phi)= \Dee(\Phi)/\tau^{-1}[m]$, which we will call the
$m$-cluster category,  
plays a similar role to
the cluster category but with respect to the combinatorics of $m$-clusters.  
We define
a bijection $W:\mphi \rightarrow \ind\Cee_m(\Phi)$.
We define an $m$-{\it tilting set} in $\Cee_m(\Phi)$ to be
a maximal set of indecomposables $\ss$ satisfying $\Ext^i_\cm(X,Y)=0$ for 
all $X,Y\in \ss$ and $i=1\dots m$.  
Then we show: 

\proclaim{Theorem 1} The map $W$ induces a bijection from 
$m$-clusters of $\Phi$ to $m$-tilting sets of $\Cee_m(\Phi)$.
\endproclaim

We then prove two facts about $\Cee_m(\Phi)$.  First: 

\proclaim{Theorem 2} The $m$-tilting sets of $\Cee_m(\Phi)$ have cardinality
$n$.  \endproclaim

Second, 
we make the natural definition of an {\it almost complete} $m$-tilting set, namely,
that it is a collection $\st$ of $n-1$ indecomposables of $\Cee_m(\Phi)$
such that $\Ext^i(X,Y)=0$ for all $X,Y\in \st$ and $i=1\dots m$,
and
we show:

\proclaim{Theorem 3} An almost complete $m$-tilting set $\st$ has $m+1$ complements.  
\endproclaim

Via Theorem 1, Theorems 2 and 3 are equivalent to facts proved about the 
$m$-cluster complex by Fomin and Reading 
([FR, Theorem 2.9 and Proposition 2.10]).  The proofs of these 
results in [FR]
depend on case-by-case
arguments and a computer check for the exceptional types.  Using Theorem 1,
our proofs of Theorems 2 and 3 therefore provide a 
type-uniform and computer-free
proof of these results from [FR] in the simply laced cases.

After this paper was completed, we received a copy of Zhu's paper [Z] which 
proves Theorems 1 and 2 without the simply laced assumption, by drawing on
some sophisticated results reported in Iyama [I], and presented 
in detail in the 
recent preprint of Iyama and Yoshino [IY].  
Our approach is different and
more elementary.


\head Clusters \endhead

We begin with a quick introduction to the combinatorics of clusters.  
Our presentation is based on [FZ] and [FR].

Let $\Phi$ be a crystallographic root system of rank $n$.  
(In fact, the assumption that $\Phi$ is crystallographic
is not essential [FR], but since we will shortly be assuming that $\Phi$ is
not merely crystallographic but also simply laced, there is no advantage
to considering the slightly more general situation.)  For convenience, we 
will also assume that $\Phi$ is irreducible; the analysis extends easily
to the reducible case.    

Label the 
vertices of the Dynkin diagram for $\Phi$ by the numbers from 1 to $n$.  
Let $W$ be the Weyl group corresponding to $\Phi$.  
Let the simple roots of $\Phi$ be $\Pi=\{\alpha_1,\dots,\alpha_n\}$, and 
let $s_i$ be the reflection in $W$ corresponding to $\alpha_i$.  

The ground set for the {\it cluster complex} $\Delta(\Phi)$ is the set 
of {\it almost positive roots}, $\aphi$, which are, by definition, the positive
roots $\Phi^+$ together with the negative simple roots $-\Pi$.  

Since the Dynkin diagram for $\Phi$ is a tree, it is a bipartite graph.
Let $I_+, I_-$ be a decomposition of $[n]$ corresponding to the bipartition.
($I_+$ and $I_-$ are determined up to interchanging $+$ and $-$.  
We fix this choice once and for all.)  

For $\epsilon \in \{+,-\}$, define the bijection 
$\tau_\epsilon:\aphi \rightarrow \aphi$ by 
$$\tau_\epsilon(\beta)=\left\{ \matrix \beta & \text { if $\beta=-\alpha_i$ for
some $i\in I_{-\epsilon}$} \\
\left(\prod_{i\in I_{\epsilon}}s_i\right)\beta &\text{ otherwise}
\endmatrix \right.$$

Now set $R=\tau_+\tau_-$.  $R$ is in some sense a deformation of a 
Coxeter element of $W$.  (We will give a 
more representation-theoretic interpretation for $R$ in Lemma 1 below.)

The crucial fact about $R$ is that every root in $\aphi$ has at least one
negative simple root in
its $R$-orbit. For that reason, the following suffices
to define a relation called compatibility.  

(c1) $-\alpha_i$ is compatible with $\beta$ iff $\alpha_i$ does not
appear when we write $\beta$ as a sum of simple roots.  (This is called the
{\it simple root expansion} for $\beta$.)

(c2) $\alpha$ and $\beta$ are compatible iff $R(\alpha)$ and $R(\beta)$ 
are compatible.  

This relation is well-defined (not {\it a priori} obvious, since
a root may have two negative simples in its $R$-orbit) and symmetric [FZ,
\S\S 3.1-2].

In fact, there is more information associated to a pair of almost  
positive roots than mere compatibility or incompatibility.  
The compatibility degree
$(\alpha \parallel \beta)$ can be defined by saying that:

(d1) $(\beta\parallel-\alpha_i)$ is the coefficient of $\alpha_i$ in the
expansion of $\beta$ if $\beta$ is positive and $0$ if $\beta$ is negative.

(d2) $(R(\beta)\parallel R(\alpha))=(\beta\parallel\alpha)$.

Compatibility degree is well-defined, and, 
if $\Phi$ is simply laced, it is also symmetric [FZ, \S 3.1].
Two roots are compatible iff their compatibility degree is zero. 

The {\it cluster complex}
 $\Delta(\Phi)$ is defined to be the simplicial complex 
whose faces are the sets of almost positive roots which
are pairwise compatible.  The facets (maximal faces) 
of $\Delta(\Phi)$ are all of
the same cardinality, $n$, the rank of $\Phi$.  They are called 
{\it clusters}.  

\head Derived Category \endhead

Fix $Q$ an orientation of the Dynkin diagram for $\Phi$.   
The representations of $Q$ are denoted $\el(Q)$.  
The bounded derived category $\Dee(Q)$ is a triangulated category, and it 
comes with a $\mathbb{Z}$ grading and 
a shift functor $[1]$ which takes $\Dee(Q)_i$ to $\Dee(Q)_{i-1}$. 
$\Dee(Q)_i$ is just a copy of $\el(Q)$. We refer to this grading as
the {\it coarse grading} on $\Dee(Q)$, and denote the degree function
with respect to this grading by $d_C$.


To give a more concrete description of $\Dee(Q)$, we 
we will define an infinite quiver $\zq$.
Its vertex set consists of 
$[n] \times \mathbb{Z}$.  For each edge from $v_i$ to $v_j$ in $Q$, 
$\zq$ has an edge from $(j,p)$ to $(i,p)$ and one
from $(i,p)$ to $(j,p-1)$, for all $p\in \mathbb{Z}$.  
This means that one way of thinking of $\zq$ is as $\mathbb{Z}$ many
copies of $\qo$ ($Q$ with its orientation reversed)
with some edges added connecting copy $i$ to copy
$i-1$.  

It turns out that $\zq$ is the Auslander-Reiten quiver for $\Dee(Q)$,
so in particular 
the indecomposables 
of $\Dee(Q)$ can be identified with the vertices of $\zq$.  See
[H] for further details and proofs.  


If $Q$ and $Q'$ are two different orientations of the 
Dynkin diagram for $\Phi$, then $\Dee(Q)$ and $\Dee(Q')$ are isomorphic
as triangulated categories, but their coarse gradings disagree.  We will 
generally therefore forget the orientation (and the grading it
induces), and write $\Dee(\Phi)$.  

Since $\Dee(\Phi)$ does not depend on the choice of an orientation, we may fix
a convenient orientation if we like.  Let $\qb$ denote the 
bipartite orientation of the Dynkin diagram of $\Phi$ in which the arrows
go from roots in $I_+$ and towards roots in $I_-$.  We want to fix a 
grading on the vertices of $\zqb$, 
which we shall call the {\it fine grading}, and 
denote it $d_F$.  
Vertices in $\zqb$ are indexed by $(i,k)$ with $i\in [n]$ 
and 
$k\in\mathbb{Z}$.  We say that a vertex $(i,k)$ is in fine degree $2k$ if
$i\in I_-$ and $2k-1$ if $i\in I_+$.  It follows that all the 
arrows in $\zqb$ diminish fine degree by 1.  

The coarse and fine gradings of $\Dee(\qb)$ are related: $d_C(M)=
\lceil d_F(M)/h \rceil$, where $h$ is the Coxeter number for 
$\Phi$.  (The Coxeter number is the order
of any Coxeter element, and can be computed from the fact that 
$|\Phi|=nh$.)  

The copy of the indecomposables of 
$\el(\qb)$ which sits in coarse degree 0, consists of the vertices of
$\zqb$ in fine degree between 0 and $-h+1$.  The indecomposables of 
$\el(\qb)$ which are projective are exactly those in fine degree 0 and
$-1$. 

Here is an example for $A_3$.  Here $h=4$, and 
$I_+$ consists of the two outside
nodes while $I_-$ is the middle node.

$$\epsfbox{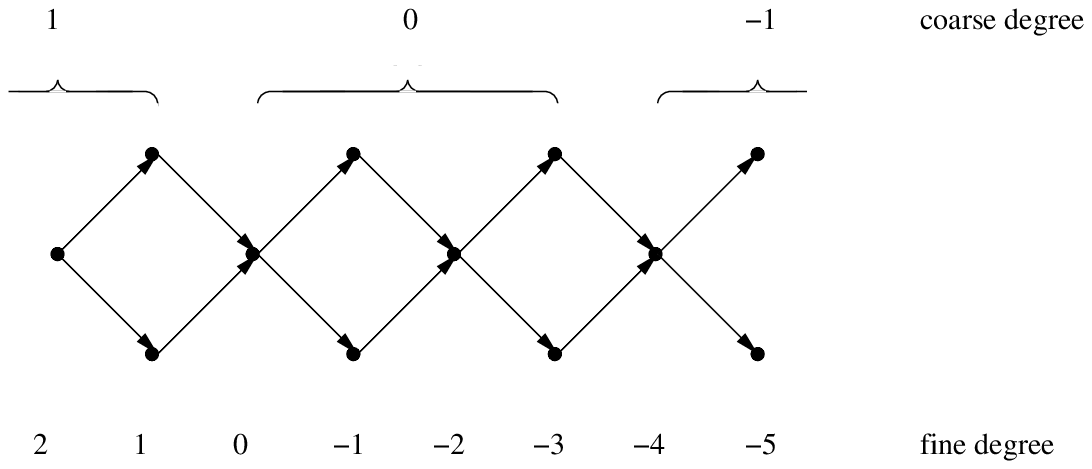}$$ 

We can define an automorphism $\tau$ of $\zq$ which  
takes 
$(i,p)$ to $(i,p+1)$ for all $i\in \{1,\dots,n\}$ and $p\in \mathbb{Z}$.
This automorphism corresponds to an auto-equivalence of $\Dee(\Phi)$,
also denoted $\tau$,  
which is the Auslander-Reiten translate for $\Dee(\Phi)$.  The functor 
$\tau$ respects the
fine degree, increasing it by 2.  The shift functor 
[1] also respects the fine degree,
decreasing it by $h$.

\head Factor categories of the derived category \endhead

Let $\Dee(H)$ be the bounded derived category of modules over a hereditary
algebra $H$, finite dimensional over a field $K$.  We quote some
general results from [BM+] and [K] about the factor of $\Dee(H)$ by a 
suitable automorphism.  

Let $G$ be an automorphism of $\Dee(H)$, satisfying conditions (g1) and (g2) 
of [BM+]:

(g1): For each $U \in \ind \Dee(H)$, only a finite number of 
$G^nU$ lie in $\ind H$ for  $n\in \mathbb{Z}$.

(g2): There is some $N\in\mathbb{N}$ such that $\{U[n]\mid U\in \ind H,
n\in [-N,N]\}$ contains a system of representatives of the orbits of $G$ on
$\Dee(H)$.  

$\Dee(H)/G$ denotes the corresponding factor category: the objects of $\Dee(H)/G$
are by definition $G$-orbits in $\Dee(H)$.  Let $X$ and $Y$ be objects of 
$\Dee(H)$, and let $\tilde X$ and $\tilde Y$ be the corresponding
objects in $\Dee(H)/G$.  Then the morphisms in $\Dee(H)/G$ are given by:
$$\Hom_{\Dee(H)/G}(\tilde X, \tilde Y)=\disj_{i\in \mathbb{Z}}\Hom_{\Dee(H)}(G^iX,Y).$$
From [K] we know that $\Dee(H)/G$ is a triangulated category, and the 
canonical map from $\Dee(H)$ to $\Dee(H)/G$ is a triangle functor.  It is shown
in [BM+, Proposition 1.2] that $\Dee(H)/G$ is a Krull-Schmidt category.  

$G$ defines an automorphism $\phi$ of the AR quiver $\Gamma(\Dee(H))$;
the factor category $\Dee(H)/G$ has almost split triangles, and the AR quiver $\Gamma(\Dee(H)/G)$ is 
$\Gamma(\Dee(H))/\phi$ [BM+, Proposition 1.3]

The shift [1] on $\Dee(H)$ passes to $\Dee(H)/G$; we use the same notation.  
Define $\Ext^i_{\Dee(H)/G}(\tilde X,\tilde Y)= \Hom_{\Dee(H)/G}(\tilde X, \tilde Y[i])$.

It is also straightforward to show [BM+, Proposition 1.4] that Serre
duality in $\Dee(H)$ passes to $\Dee(H)/G$, so  
$\Ext^1_{\Dee(H)/G}(\tilde X,\tilde Y)$ is dual to 
$\Hom_{\Dee(H)/G}(\tilde Y,\tau \tilde X).$

\head Cluster category \endhead

The cluster category is defined by $\Cee(H)=\Dee(H)/\tau^{-1}[1]$.
Since $\tau^{-1}[1]$ is an automorphism satisfying conditions (g1) and
(g2), $\Cee(H)$ is a triangulated Krull-Schmidt category. 

We will mainly be concerned in this section 
with $\Cee(\Phi)=\Dee(\Phi)/\tau^{-1}[1]$.
In [BM+], the connection between $\Cee(\Phi)$ and $\Delta(\Phi)$ is made 
via a reformulation of clusters in terms of {\it decorated representations}
due to Marsh, Reineke, and Zelevinsky [MRZ].
We proceed in a different fashion, the
basis of which is our Lemma 1 below.  

Fix $Q=\qb$.  
Identify the indecomposables of $\Dee(\Phi)$ with the vertices of $\zq$.  
It is clear that the
vertices of 
$\zq$ satisfying $2\geq d_F(M) \geq -h+1$ are a fundamental domain for
$\tau^{-1}[1]$.  The representations of $Q$ in coarse degree $0$ correspond
to indecomposables with fine degree between 0 and $-h+1$, so the 
fundamental domain we have identified 
for $\tau^{-1}[1]$ consists of the indecomposable 
representations of $Q$ in coarse degree zero 
together with $n$ extra indecomposables 
which correspond to the 
injective representations in coarse degree $1$.

We wish to put these indecomposables in bijection with $\aphi$. 
Given a representation $V$ of $Q$, its dimension is by definition 
$\dim(V)=\sum_i \dim_K (V_i) \alpha_i$.  
By Gabriel's Theorem, $\dim$ is a bijection from indecomposable 
representations of $Q$ to $\Phi^+$.  
We write $V(\beta)$ for the indecomposable
representation in coarse degree zero 
whose dimension is $\beta$.  We write $P_i$ for the 
projective representation corresponding to vertex $v_i$, and we write
$I_i$ for the injective representation correspoding to vertex 
$v_i$.  Observe that $\tau P_i=I_i[-1]$.  
We define $V(-\alpha_i)$ to be $I_i[-1]$.  

\proclaim{Lemma 1} $V(R(\alpha))=V(\alpha)[1]$. \endproclaim

\demo{Proof} 
On the representations of $Q$ which do not lie in fine degree
0 or $-1$ (i.e. the indecomposable 
representations of $Q$ which are not projective), $\tau_+\tau_-$ acts like a product of the corresponding reflection
functors, and the product of the reflection functors coincides with
$\tau$ [BB].
So $V(\tau_+\tau_-(\beta))=\tau(V(\beta))=V(\beta)[1]$, as desired.  

Now consider the case that $V(\alpha)$ is projective.  If $V(\alpha)=P_i$
is simple projective, then $i\in I_-$ and $\alpha=\alpha_i$, so
$\tau_+\tau_-(\alpha)=-\alpha_i$.  If $V(\alpha)=P_i$ is non-simple projective,
then $i\in I_+$ and $\alpha=\alpha_i+$(the sum of the adjacent roots).  Thus,
again, $\tau_+\tau_-(\alpha)=-\alpha_i$.  In both these cases, $V(R(\alpha))=
V(-\alpha_i)=I_i[-1]=\tau P_i=\tau V(\alpha)=V(\alpha)[1]$, as desired.

Finally we consider the case where $\alpha=-\alpha_i$.  
For $i\in I_-$, we know that $V(-\alpha_i)$ sits in fine degree $2$.
Now $\tau_+\tau_-(-\alpha_i)=\tau_+(\alpha_i)=\alpha_i+$(the sum of the
roots adjacent to $\alpha_i$).  We recognize this as $\dim I_i$: 
in other words, $V({R(\alpha)})=I_i=V(\alpha)[1]$, as
desired.  
For $i\in I_+$, the object $V(-\alpha_i)$ sits in fine degree 1.  In this case
$\tau_+\tau_-(-\alpha_i)=\alpha_i$.  Now $V(\alpha_i)=I_i$, 
so again $V({R(\alpha)})=\tau V(\alpha)=
V(\alpha)[1]$.    This completes the proof.  
\enddemo

The connection between representation theory and clusters now appears
strongly:

\proclaim{Proposition 1 [BM+]} 
$\dim_K\Ext^1_{\Cee(\Phi)}(V(\beta),V(\alpha))=(\beta\parallel\alpha)$ \endproclaim

\demo{Proof} We check the two defining properties of compatibility degree
given above. 
$$\align
\text{(d1)\quad}\dim_K\Ext^1_{\Cee(\Phi)}(V(\beta),V(-\alpha_i))&=
\dim_K\Ext^1_{\Cee(\Phi)}(V(\beta),I_i[-1])\\
&=\dim_K\Hom_{\Cee(\Phi)}(I_i[-1],\tau V(\beta))\\
&=\dim_K\Hom_{\Cee(\Phi)}(\tau^{-1}I_i[-1],V(\beta))\\
&=\dim_K\Hom_{\el(\Phi)}(P_i,V(\beta))
\endalign$$ 
(The first equality is because $V(-\alpha_i)=I_i[-1]$.  
The second equality is by Serre duality.  The third follows because $\tau$
is an auto-equivalence; the fourth from the fact that $\tau^{-1}I_i[-1]=P_i$.)
Now
$\dim_K\Hom_{\el(Q)}(P_i,V(\beta))$ is the coefficient
of $\alpha_i$ in the simple root expansion of $\beta$, proving condition (i).

(d2) The invariance under $R$ follows by Lemma 1 
from fact that [1] is an auto-equivalence 
of $\Cee(\Phi)$.
\enddemo

Thus, the roots in a cluster correspond to a maximal collection of
irreducible modules in $\Cee(\Phi)$ such that all the $\Ext_{\Cee(\Phi)}^1$'s between them
vanish.  This is exactly the definition of a {\it (cluster-)tilting set} for 
$\Cee$, so we have seen that {tilting sets} for $\Cee$ are in
one-one correspondence with clusters for $\Phi$.

\head $m$-Clusters \endhead

The $m$-clusters are a simplicial complex whose ground set, denoted 
$\mphi$, consists of the negative simple roots $-\Pi$ together with $m$
copies of $\Phi^+$.  These $m$ copies are referred to as having 
$m$ different
``colours'' $1$ through $m$.  
To keep track of the roots of different colour, we use
superscripts.  So $\beta^i$ is the root $\beta$ with colour $i$.  
Negative simple roots are considered to have 
colour 1.  

Fomin and Reading define an $m$-ified rotation on $\mphi$:

$$R_m(\alpha^k)=\left\{ \matrix \alpha^{k+1} &\text{ if $\alpha\in\Phi^+$ 
and $k<m$} \\
R(\alpha)^1 &\text{ otherwise}\endmatrix \right.$$

Again, the crucial fact (which follows from the $m=1$ case) 
is that every root has at least one
negative simple in its $R_m$-orbit.  

We now follow [FR] in defining a relation called compatibility.  
(Strictly speaking, perhaps,
we should call this $m$-compatibility, but no ambiguity will result,
because this is a relation on $\mphi$, not $\aphi$.)

(m1) $-\alpha_i$ is compatible with all negative simple roots and any 
positive root (of whatever colour) 
that does not use $\alpha_i$ in its simple
root expansion.

(m2) $\alpha$ and $\beta$ are compatible iff $R_m(\alpha)$ and $R_m(\beta)$
are compatible.  

Because of the crucial fact mentioned above, this is sufficient to define
compatibility, but not to prove that such a relation exists.  This is 
verified in [FR], where the relation 
is also shown to be symmetric.  
We shall give our own proof that there is a compatibility relation on
$\mphi$ satisfying (m1) and (m2), below (Proposition 2).  

The $m$-cluster complex $\Delta_m(\Phi)$ is the simplicial complex on
$\mphi$
whose
faces are sets of pairwise compatible roots.  The facets of the
complex are called $m$-clusters.  



\head $m$-Cluster Category \endhead

We define the $m$-cluster category to be 
$\Cee_m(\Phi)=\Dee(\Phi)/\tau^{-1}[m]$.  
This category is discussed in [K, Section 8.3], where it is shown
to be triangulated, and in [KR].  
It is also being studied at present by
A. Wraalsen [W].  
The type $A$ case has been considered in detail in [BM].  


We now identify the indecomposables of $\Cee_m(\Phi)$ with $\mphi$, as 
follows.  
For $\beta^j$ a positive root in $\mphi$, let $W(\beta^j)=V(\beta)[j-1]$.
Let $W(-\alpha_i)=I_i[-1]$.
Observe that the set of $W(\beta^k)$ which we have  
identified are a fundamental domain
with respect to $F=\tau^{-1}[m]$, and therefore they correspond in a 
1-1 fashion to the indecomposables of $\Cee_m(\Phi)$.  


We now prove the $m$-ified analogue of Lemma 1.

\proclaim{Lemma 2} $W (R_m(\beta^k))=W(\beta^k)[1]$. \endproclaim

\demo{Proof} 
There are three cases to consider: firstly when $\beta^k=-\alpha_i$, 
secondly when $\beta$ is a positive
root and $k<m$, and thirdly when $\beta$ is a positive root and $k=m$.  

In the first case, $\beta^k=-\alpha_i$, and $R_m(-\alpha_i)=
R(-\alpha_i)^1$.  In this case $W(-\alpha_i)=I_i[-1]$,
and by the proof of Lemma 1, $W(R(-\alpha_i)^1)=I_i$, which proves
the claim in this case.  

In the second case ($\beta$ a positive root and $k<m$), $R_m(\beta^k)=
\beta^{k+1}$, and the desired result follows by the definition of
$W$.  

In the third case, $R_m(\beta^m)=R(\beta)^1$.  By the proof of Lemma 1, 
$W(R_m(\beta^m))=\tau(V(\beta))=V(\beta)[m]=W(\beta^m)[1]$, as desired.  
\enddemo


We now prove the following analogue of Proposition 1.7(b) of [BM+]:

\proclaim{Lemma 3} If $X,Y\in \Cee_m(\Phi)$, then $\dim_K\Ext^i
_{\ceem}(X,Y)=
\dim_K\Ext^{m+1-i}_{\ceem}(Y,X)$. \endproclaim

\demo{Proof}
This is essentially (a slightly naive version of) the Calabi-Yau condition of
dimension $m+1$, proved for $\Cee_m(\Phi)$ by Keller in [K, Section 8.3].  
Observe that $\dim_K\Ext^{i}_\ceem(X,Y)=\dim_K\Ext^1_\ceem(X,Y[i-1])=
\dim_K\Hom_\ceem(\tau^{-1} Y[i-1],X)
=\dim_K\Hom_\ceem(Y[i-1-m],X)=\dim_K\Ext^{m+1-i}_\ceem(Y,X).$ The second equality follows
by Serre duality, and the third because we are in $\Dee(\Phi)/\tau^{-1}[m]$.
\enddemo

We now prove an $m$-ified analogue of Proposition 1.  Here, we consider
only compatibility, not compatibility degree, as [FR] does
not define compatibility degree in the $m$-cluster context (though it
is not difficult to do so).  

\proclaim{Proposition 2} There is a well-defined and symmetric 
relation (called compatibility)
on $\mphi$ satisfying (m1) and (m2).
A pair of coloured roots
$\beta^j$ and $\gamma^k$ are compatible in
$\mphi$ iff $\Ext_{\Cee_m(\Phi)}^i(W(\beta^j),W(\gamma^k))=0$ for $i=1\dots m$.  \endproclaim

\demo{Proof} The proof is analogous to the proof of Proposition 1:
since (m1) and (m2) define a unique relation, it suffices to check
that the relation given in the statement of the proposition 
is symmetric and satisfies (m1) and (m2).  Symmetry follows from Lemma 3.  
We now proceed to check (m1) and (m2).  

$$
\align
\text{(m1)\quad}\dim_K\Ext^j_{\Cee_m(\Phi)}(W(\beta^k),W(-\alpha_i))=&\dim_K\Ext^j_{\Cee_m(\Phi)}
(W(\beta^k),I_i[-1])\\=&\dim_K\Ext^1_{\Cee_m(\Phi)}(W(\beta^k),I_i[j-2]) \\
    =&\dim_K\Hom_{\Cee_m(\Phi)}(I_i[j-2],\tau W(\beta^k))\\
=&\dim_K\Hom_{\Cee_m(\Phi)}(\tau^{-1}I_i[j-2],W(\beta^k))\\
=&\dim_K\Hom_{\Cee_m(\Phi)}(P_i[j-1],W(\beta^k)).
\endalign$$
(The second equality is from the definition of $\Ext^i$; the third is from
Serre duality; the fourth follows because $\tau$ is an auto-equivalence.)
If $k\ne j$, then $\Hom_{\Cee_m(\Phi)}(P_i[j-1],W(\beta^k))$
 is zero.  If $k=j$, it is the coefficient of 
$\alpha_i$ in the root expansion of $\beta$.
This implies that $\beta^k$ is incompatible with $-\alpha_i$ if and only if
$\alpha_i$ appears with positive coefficient in the simple root expansion
of $\beta$.  

(m2) This follows from Lemma 2 and the fact that $[1]$ is an auto-equivalence
of $\Cee_m(\Phi)$.  
\enddemo

We define an $m$-{\it tilting set} in $\Cee_m(\Phi)$ to be
a maximal set of indecomposables $\ss$ satisfying $\Ext^i(X,Y)=0$ for 
all $X,Y\in \ss$ and $i=1\dots m$.  
The following theorem is an immediate consequence of Proposition 2:

\proclaim{Theorem 1} The map $W$ induces a bijection from 
$m$-clusters of $\Phi$ to $m$-tilting sets of $\Cee(\Phi)$.
\endproclaim

\head Combinatorics of $m$-clusters \endhead

In this section, we prove Theorems 2 and 3.  These are
representation-theoretic
reformulations
of the following results from [FR]:

\proclaim{Theorem $2'$ [FR, Theorem 2.9]} All the facets of 
$\Delta_m(\Phi)$ are of
size $n$.
\endproclaim

\proclaim{Theorem $3'$ [FR, Proposition 2.10]} Given a set $\st$ of $n-1$ pairwise compatible 
roots from $\mphi$, there are exactly $m+1$ roots not in $\st$ which 
are compatible with all the roots of $\st$.  
(In other words, every 
codimension 1 face of $\Delta_m(\Phi)$ is contained in exactly $m+1$
facets.)
\endproclaim

The proofs of these results in [FR] rely on the following theorem,
which is proved on a type-by-type basis, with a computer check for
the exceptionals.  
We will give a type-free proof.

\proclaim{Theorem 4 [FR, Theorem 2.7]} If $\alpha$ and $\beta$ 
are roots of $\Phi$ contained in a parabolic root system $\Psi$ 
within $\Phi$, then
$\alpha^{(i)}$ and $\beta^{(j)}$ are compatible in $\mphi$ iff they are
compatible in $\mpsi$.  \endproclaim

\demo{Proof} Let $\bar X$ and $\bar Y$ be indecomposables of $\ceem$ corresponding
to $\alpha^{(i)}$ and $\beta^{(j)}$ respectively.  
Let 
$X$ and $Y$ be corresponding indecomposables of $\dd$, chosen so that
$2 \geq d_F(X), d_F(Y) \geq -mh+1$.  Without loss of generality, we may assume that 
$d_F(X)\geq d_F(Y)$.  

Suppose that 
$\alpha^{(i)}$ and $\beta^{(j)}$ are not compatible in $\mphi$.  So there
is some $1\leq k \leq m$ with $\Ext_\ceem^k(\bar X,\bar Y)\ne 0$.
This asserts that there is a non-zero 
morphism in $\dd$ from some $G^pX$ to $Y[k]$,
where $G=\tau^{-1}[m]$ and
$p$ is some integer.  Since
$d_F(X)\geq d_F(Y)$, so $d_F(X)-d_F(Y[k])\geq h$, it follows that
 $p$ must be strictly positive.  On the other hand, $d_F(Y[k])>d_F(G^pX)$
for $p\geq 2$.  So $p=1$, and $\Hom_\dd(\tau^{-1}X[m],Y[k])\ne 0$.  
By Serre duality, $\Ext^1_\dd(Y[k],X[m])\ne 0$, so $\Ext^{m+1-k}_\dd(Y,X)
\ne 0$.  The crucial point here is that we know this statement
on the level of the derived category, rather than just the 
$m$-cluster category.  

Let $Q'$ be the subquiver of $Q$ corresponding to $\Psi$.  
There is a 
natural inclusion of $\el( Q')$ into $\el(Q)$ as a full subcategory, 
which extends to an inclusion of $\Dee(Q')$ into $\Dee(Q)$ as a full
triangulated subcategory, where the inclusion respects the coarse grading.
$X$ and $Y$ represent $\alpha^{(i)}$ and $\beta^{(j)}$ respectively
in both $\Cee_m(\Phi)$ and $\Cee_m(\Psi)$.  
Thus, the non-vanishing $\Ext$ that we have shown exists in $\Dee(\Phi)$ 
also exists in $\dpsi$, and testifies that
$\alpha^{(i)}$ and $\beta^{(j)}$ are not compatible in $\mpsi$ either.  

The converse is proved similarly.  
\enddemo

Now that we have established Theorem 4, the proofs of Theorems
2 and 3 
go through exactly as in [FR].  We 
include the proofs for completeness.

\demo{Proof of Theorem 2} The proof is by induction
on $n$.  The statement is clear when $n=1$.  
Let $\ss$ be an $m$-tilting set in $\Cee_m(\Phi)$. 
Pick $X$ an indecomposable in $\ss$.  Applying $\tau$ if necessary,
we may assume that $X$ is of the form $I_i[-1]$ for some $i$.  Let
$Q'$ be the quiver $Q$ with the vertex $v_i$ removed, and let
$\Psi$ be the associated root subsystem.  
For each indecomposable $Y \in \ss\setminus \{X\}$, choose a representative
$\hat Y$
in $\Dee(Q)$ with fine degree between 2 and $-hm+1$ (in other words,
$\hat Y$ is either of the form $I_j[-1]$ for $j\ne i$ 
or in $\el(Q)[k]$ for some $0\leq k
\leq m-1$.  

Since $Y$ is compatible with $X$, $\Ext^j_\ceem(Y,I_i[-1])=0$ for all $1\leq j \leq m$.
By Serre duality, this is equivalent to the condition that
$\Hom_\ceem(P_i[j-1],Y)=0$ for all $1\leq j \leq m$, or, in other words,
that, if $\hat Y\in \el(Q)[k]$, that in fact $\hat Y\in \el(Q')[k]$.
Thus, 
by Theorem 4, the images of the $\hat Y$ form
an $m$-tilting set in $\Cee_m(\Psi)$, so $\ss\setminus \{X\}$ contains $n-1$
indecomposables by induction, and thus $\ss$ contains $n$ indecomposables.   
\enddemo

\demo{Proof of Theorem 3} The proof is, again, by induction on
$n$.  The base case, when $n=1$, is clear.  For the induction step,
let $\st$ be an almost complete 
tilting set.  As before, we choose an indecomposable 
 $X$ in $\st$, which we may assume is of the form $I_i[-1]$, and then
we observe that $\st\setminus\{X\}$ consists of an almost complete tilting set
for a root system of rank $n-1$, and the $m+1$ complements for that
almost complete $m$-tilting set 
are precisely the complements of $\st$ in
$\Cee_m(\Phi)$.  
\enddemo

\head Acknowledgements \endhead

We would like to thank Colin Ingalls for sharing his knowledge
of quiver representations, and  Drew Armstrong, 
Sergei Fomin, Osamu Iyama, Bernhard Keller, Nathan Reading, Idun Reiten,
Ralf Schiffler, David Speyer,
Andrei Zelevinsky, and an anonymous referee for helpful comments.

\enddocument